\documentclass {amsart}
\usepackage{amsmath,amsthm,amsfonts,amssymb}
\usepackage{a4}
\newtheorem{thrm}{Theorem}[section]
\newtheorem*{thrm*}{Theorem}
\newtheorem*{kor*}{Korollar}
\newtheorem*{cor*}{Corollary}
\newtheorem*{prop*}{Proposition}
\newtheorem{defin}[thrm]{Definition}
\newtheorem{prop}[thrm]{Proposition}
\newtheorem{lemma}[thrm]{Lemma}

\newtheorem{remark}[thrm]{Remark}

\def\Ad{{\rm Ad}}
\def\ad{{\rm ad}}
\def\SO{{\rm SO}}

\def\Sp{{\rm Sp}}
\def\U{{\rm U}}
\def\GL{{\rm GL}}

\def\conv{{\rm conv}}
\def\emb{{\rm emb}}
\def\cohom{{\rm cohom}}
\def\Hom{{\rm Hom}}
\def\conv{{\rm conv}}

\def\rk{{\rm rk}}
\def\iso{{\rm Iso}}
\def\exp{{\rm exp}}
\def\End{{\rm End}}
\def\kdot{\!\cdot\!}
\def\T{T^\mathbb{C}}
\def\G{G^\mathbb{C}}
\def\t{\frak t^\mathbb{C}}
\def\P{\mathbb{P}}
\def\g{\frak g^\mathbb{C}}
\begin{document}
\title{Polar Actions on Hermitian and  Quaternion-K\"ahler Symmetric  Spaces } 
\author{Samuel Tebege}
\address{Samuel Tebege, Mathematisches Institut, Universit\"at zu K\"oln, Weyertal 86-90, 50931 K\"oln, Germany}
\email{stebege@math.uni-koeln.de}
\subjclass{53C55, 53C26, 57S15}
\date{January 17, 2007}
\keywords{Polar actions, Hermitian and quaternion-K\"ahler symmetric spaces}
\begin{abstract}

We analyze polar actions on Hermitian and quaternion-K\"ahler symmetric spaces of compact type. For complex integrable polar actions on Hermitian symmetric spaces of compact type we prove a reduction theorem and several corollaries concerning the geometry of these actions. The results are independent of the classification of polar actions on Hermitian symmetric spaces. In the second part we prove that polar actions on Wolf spaces are quaternion-coisotropic and that isometric actions on these spaces admit an orbit of special type, analogous to the existence of a complex orbit for an isometric action on a compact homogeneous simply-connected K\"ahler manifold.
\end{abstract}
\maketitle
\section{Introduction}
Let $M$ be a complete Riemannian manifold and $G$ a compact connected  Lie group acting isometrically on $M$. A submanifold $\Sigma \subset M$ is called a \emph{section} if $\Sigma$ intersects
every $G$-orbit and is orthogonal to the $G$-orbits in all common
points. An isometric $G$-action that admits a section is called \emph{polar}. A
polar action, such that the section $\Sigma$ is flat, is called \emph{hyperpolar}.

The first systematic study of polar actions on Riemannian manifolds began with the paper {\it A general theory of canonical forms}  of Palais and Terng \cite{Palais1}. Nevertheless a special case, namely the hyperpolar ones, were considered a long time before. For example Conlon has shown that hyperpolar actions (he called them actions admitting a $K$-transversal domain) are variationally complete \cite{Conlon}, generalizing results of Bott and Samelson, who have proven that the adjoint action of a compact Lie group and the isotropy action of a compact symmetric space are variationally complete,  these two actions being particularly simple examples of polar actions \cite{Bott}. 

The case of polar representations has been analyzed by Dadok, who has classified them thereby proving that they are orbit equivalent to $s$-representations \cite{Dadok}. Since the slice representation of an arbitrary polar action is again polar, this shows that the geometry of polar actions is intimately connected with the geometry of symmetric spaces. Actually most of the known examples of polar actions are given on symmetric spaces. Now all these examples on compact irreducible symmetric spaces of higher rank turned out to be even hyperpolar. The question whether this holds in full generality remains open, however the question has been answered for a large class of spaces.  The first answer to this question was given on certain subclasses of Hermitian symmetric spaces. Podest\`a and Thorbergsson have classified polar actions on the complex quadric making use of their result that in this case the action is coisotropic, meaning that the sections are totally real submanifolds \cite{Podesta}. The same idea was used by Gori and Biliotti for the complex Grassmannian \cite{Gori} and by Biliotti for the remaining Hermitian symmetric spaces of higher rank \cite{Biliotti}. It turned out that polar actions on irreducible higher rank Hermitian symmetric spaces of compact type are hyperpolar.  

Recently Kollross finished the classification of polar actions on symmetric spaces of type I showing that polar actions are hyperpolar on these spaces if they are of higher rank. In principle his methods should also work for symmetric spaces of type II but the computational complexity becomes quite large, so this case remains open. 

The first part of this paper is concerned with polar actions on Hermitian symmetric spaces. The main motivation for us was to supplement geo\-metric ideas which might be helpful in a conceptual proof of the fact mentioned above making use of the rich geometrical structure of these spaces. In particular our results are independent of the classification results of polar actions on Hermitian symmetric spaces. Especially the existence of a momentum map for isometric actions turns out to be helpful. Let $H$ be a compact connected Lie group and recall that for a symmetric space the polarity of the $H$-action is equivalent to the integrability of the principal horizontal distribution $$M_r\to M_r/H,$$ where $M_r$ denotes the set of regular points of the $H$-action \cite{Palais2}. On a Hermitian symmetric space the result of Podest\`a and Thorbergsson mentioned above implies that the principal horizontal distribution is totally real. We call a polar action on a K\"ahler manifold \emph{complex integrable} if the complexified principal horizontal distribution is integrable. The main result of the first part is a reduction theorem. We reduce the acting group to a certain subtorus which is defined by the momentum map while reducing at the same time the codimension of the section. So let $M$ be an irreducible Hermitian symmetric space of compact type endowed with a complex integrable polar action of a compact connected Lie group $H$ and $\mu:M\to \frak h^*$ a momentum map for the $H$-action. We will prove:
\begin{thrm*}[Reduction Theorem]
Let $p\in M$ be a regular point of the $H$-action, $\Sigma$ the section through $p$ and $T:=H_{\mu(p)}/H_p$. Then:
\begin{enumerate}\renewcommand{\labelenumi}{(\roman{enumi})}\item
$\Sigma^{\mathbb{C}}:=\overline{T^{\mathbb{C}}p}$ is a totally geodesic submanifold of $M$ which contains $\Sigma$ as a Lagrangian submanifold. \item The $T$-action on $\Sigma^\mathbb{C}$ is polar with $\Sigma$ as a section. \item  The complexified section $\Sigma^\mathbb{C}$ splits as a product of complex projective spaces, i.~e.~$\Sigma^\mathbb{C}=\mathbb{P}_{n_1}(\mathbb{C})\times \cdots \times \mathbb{P}_{n_k}(\mathbb{C})$.\end{enumerate}
\end{thrm*}
In the next section we will give applications of this theorem. We obtain several corollaries. The first one is a criterion for the hyperpolarity of the $H$-action.
\begin{cor*}The $H$-action is hyperpolar if and only if there exists a totally geodesic  momentum fiber through a regular point of the $H$-action. 

\end{cor*}
Moreover the Reduction Theorem enables us to  show that the sections of the $H$-action  are necessarily products of real projective spaces. Compare this also to Theorem 5.4 of \cite{Kollross2}, where Kollross proves that a section of a polar action is covered by a product of spaces of constant curvature.
\begin{cor*}Let $\Sigma$ be a section for the $H$-action on $M$. Then $\Sigma$ is isometric to a product of real projective spaces, i.~e.~$\Sigma=\mathbb{P}_{n_1}(\mathbb{R})\times \cdots \times \mathbb{P}_{n_k}(\mathbb{R})$. The $H$-action is hyperpolar if and only if $n_1=\cdots=n_k=1$.
\end{cor*}
Observing that every polar action on a complex  projective space is complex integrable this leads to a new conceptual proof of the fact that in this case a section must be a real projective space.  This was first observed by Podest\`a and Thorbergsson in \cite{Podesta2}, where they classified polar actions on the rank one spaces.
\begin{cor*} Let $X$ be a complex projective space and $H$ a connected compact Lie group acting polarly on $X$. Then the action is complex integrable and a section is isometric to a real projective space.
\end{cor*}
We also prove that complex integrable polar actions are provided by polar $\mathbb{C}$-asystatic ones.
\begin{prop*}Let $M$ be a compact homogeneous simply-connected irreducible K\"ahler manifold and $H\times M\to M$ a polar $\mathbb{C}$-asystatic action of a compact connected Lie group $H$. Then the $H$-action is complex integrable.
\end{prop*}
One of the main difficulties in a further analysis is to understand the geometry of the torus actions of $T=H_{\mu(p)}/H_p$ and $\T$. What is the relation between a regular point $p\in M$ and the torus $T=H_{\mu(p)}/H_p$, respectively $\T$? 
While the general case remains open due to fact that algebraic torus actions on generalized flag manifolds are only well understood if the torus in question is maximal, we obtain some further results for the series {\bf {\sl C I}}, i.~e.~$M=\Sp(n)/U(n)$. The key property of this spaces we will use is that they are symmetric spaces of maximal rank, i.~e.~$\rk\; G/K=\rk\; G$. First of all we identify the generic torus orbit type for $M$ using results of Dabrowski \cite{Da}.
\begin{prop*} The generic torus orbit type for $M=G/K=\Sp(n)/\U(n)$, i.~e.~the generic torus orbit type of a maximal torus $\T\subset G^\mathbb{C}$ is isomorphic to the product of $n$ complex projective lines. 
\end{prop*}
Note that this is the orbit type for the torus which is defined by the momentum map for the isotropy action on $M$, i.~e.~$T=K_{\mu(p)}/K_p$. It follows that for a hyperpolar action of maximal cohomogeneity the point $p$ is contained in the non-degenerate stratum $Z(\T)$ of the Torus $\T$. This is the stratum which is defined by the intersection of the big cells of all Borel subgroups in $G^\mathbb{C}$ which contain $\T$.
\begin{prop*}\label{Za} Let $H\times M\to M$ be a complex integrable hyperpolar action with maximal cohomogeneity of a compact connected Lie group on $M=\Sp(n)/\U(n)$. Then every regular point $p\in M$ is contained in the non-degenerate stratum $Z(\T)$ of the complexification of the torus $T=H_{\mu(p)}/H_p$.
\end{prop*}
The last section of this paper is concerned with polar actions on quaternion-K\"ahler symmetric spaces of compact type, i.~e.~the so-called Wolf spaces.
The first result generalizes the main result of \cite{Podesta} to the quaternionic case.
\begin{thrm*}\label{qco}Let $M$ be a Wolf space and $H$ a connected compact Lie group acting polarly on $M$. Then the $H$-action is quaternion-coisotropic. 
\end{thrm*}
For this result we will make essential use of the classification result of Kollross \cite{Kollross2}. The main difficulty in generalizing the proof of the analogous result for the K\"ahler case is that we cannot guarantee the existence of a quaternionic orbit. However we will show:
\begin{thrm*}Let $M$ be a Wolf space and $H$ a compact connected Lie group acting isometrically on $M$. Then the $H$-action admits at least one of the following orbits:
\renewcommand{\labelenumi}{(\roman{enumi})}
\begin{enumerate}\item a complex orbit,
\item a totally geodesic quaternionic orbit,
\item a $\mathbb{Z}_2$-quotient of a complex orbit in the twistor space $Z$ of $M$.
\end{enumerate}
\end{thrm*}
\noindent
{\it Acknowledgments.} I would like to thank my teacher Professor Gudlaugur Thorbergsson, who supported me since the very beginning of my studies. Without his patience and constant encouragement this work would not have been possible.
\section{Basic properties of polar actions}
\noindent In the following we will collect some well-known facts about polar actions. The main references are \cite{Palais1} and \cite{Palais2}.
First of all it is easy to see that for all $g\in G$ the translated section $g\Sigma$ is again a section. Together with the fact that $\Sigma$ intersects every orbit, we see that there exists a section through every point of $M$. If $p\in M$ is a regular point of the $G$-action then the section through $p$ is uniquely determined. It is just the image of the normal space of the orbit $G\kdot p$ under the exponential map of $M$. Given any polar action, the slice representation, i.~e.~the restriction of the isotropy representation to the normal space, is again polar. If $\Sigma$ is a section through $p$ then $T_p\Sigma$ is a section for the slice representation. The orbit structure of polar representations is well understood. Dadok \cite{Dadok} has proven that every polar representation of a compact connected Lie group is orbit equivalent to the isotropy representation of a symmetric space. Together with the facts about the slice representation mentioned above  this turns out to be an important tool for the analysis of polar actions. It is also well known, that a section is automatically totally geodesic. A proof can be found in 
\cite{Tebege}.
\section{Polar actions on Hermitian symmetric spaces}
\subsection{Hermitian symmetric spaces as generalized flag manifolds}
The main references for this section are \cite{Helgason}, \cite{Borel}, \cite{Humphreys}. Let $M$ be an irreducible Hermitian symmetric space of compact type. Then $M$ can be realized as a quotient $G/K$, where $G$ is a connected compact simple Lie group with trivial center and $K$ is a maximal connected proper subgroup of $G$ with discrete center. There are two other important geometric realizations of $M$, which we will use in the sequel: one as an adjoint orbit and the other as a generalized flag manifold.

Let $\frak g=\frak k\oplus\frak m$ be the corresponding Cartan decomposition of $\frak g$. Then there exists an element $z$ in the center of $\frak k$ such that $\ad_z:\frak m\to\frak m$ gives the complex structure of $M$. Since $z$ lies in the center of $\frak k$, it is easy to see that $K$ is the stabilizer in $z$ for the adjoint action of $G$ on $\frak g$ and hence $M=\Ad(G)\kdot z$. This description as an adjoint orbit will guarantee that later on all isometric actions of  Lie subgroups of $G$ on $M$ are automatically Hamiltonian.

On the other hand let $\G$ denote the complexification of $G$. Then $M$ is a $\G$-homogeneous space and hence $M=\G/P$, where $P$ is a parabolic subgroup of $\G$. 

There exists the following root-theoretical description of $P$: We fix a Cartan subalgebra $\t\subset \g$ and consider the root system $\Delta=\Delta(\g,\t)$ of $\g$ with respect to $\t$. Let $\Delta^+$ denote the set of positive roots for some ordering of $\Delta$ and $\pi\subset \Delta^+$ the set of simple roots. Then $$\frak b=\t\oplus\bigoplus_{\alpha\in\Delta^+}\g_\alpha$$
is a Borel subalgebra and the parabolic algebras $\frak p$ with $\frak p\supset\frak b$ are in one-to-one correspondence with subsets $\pi_P$ of $\pi$. The parabolic subgroups corresponding to Hermitian symmetric spaces are maximal ones and hence in this case the subset $\pi_P$ arises from $\pi$ by deleting one root.  More precisely: The sets $\pi_P=\pi\setminus\{\alpha\}$ which give rise to Hermitian symmetric spaces are characterized by the condition that the coefficient of $\alpha$ in the sum description of the highest root of $\Delta$ is one. Let $\omega_i$ denote the fundamental dominant weight corresponding to the root $\alpha_i\in\pi$ and $V$ a $\G$-module with heighest weight $$\delta=\sum_{\alpha_i\in \pi_P}\omega_i\; \text{  and heighest vector }v_\delta\in V,$$
then we have  $$\G/P\cong\G\kdot[v_\delta]\subset \mathbb{P}(V).$$
Let $B\subset\G$ be a Borel subgroup. The $B$-orbits in $M$ are called Schubert cycles. From the Bruhat-decomposition of $\G$ it follows that there exists a unique open and dense $B$-orbit in $M$, which is called the big cell of $B$ in $M$. The Weyl group $W:=N_\G(\T)/\T$ operates simply transitive on the set $\mathcal{B}^{\T}$ of all Borel subgroups containing $\T$, so we can write $$\mathcal{B}^{\T}=\{wB^-w^{-1}\mid w\in W\},$$
where $B^-$ is the opposite Borel subgroup of $B$, i.~e.~$B\cap B^-=\T.$ For each of this Borel subgroups one gets one big cell and the intersection of all these big cells is called the non-degenerate stratum $Z(\T)$. The big cell for $B^-$ is just the $B^-$-orbit through the origin $eP\in G/P$ and hence one has the following description:
$$Z(\T)=\bigcap_{w\in W}wB^-w^{-1}wP=\bigcap_{w\in W}wB^-P=\bigcap_{w\in W/W_P}wB^-P,$$ where $W_P$ is the subgroup of $W$ generated by reflections $\alpha\in \pi_P$ or alternatively $W_P=N_P(\T)/\T$.
\subsection{Torus actions and toric varieties}
The main reference for this section is \cite{Oda}. Since this is non-standard material from the differential geometric viewpoint, we allow ourselves to be a little more elaborate. The key idea is that toric varieties are in one-to-one correspondence to combinatorial objects called fans, which are families of certain cones. Let us first define the notion of a toric variety:
\begin{defin} Let $X$ be an irreducible normal scheme $X$ over $\mathbb{C}$ which is locally of finite type over $\mathbb{C}$ and seperated. Then $X$ is called a toric variety, if it containes an algebraic torus $\T$ as an open set, such that $X$ admits an algebraic action $\T\times X\to X$ which extends the group multiplication $\T\times\T\to\T$.
\end{defin}

We now would like to define the involved combinatorial objects.
Let $N$ be a free $\mathbb{Z}$-module of rank $r$, $M$ its dual module, i.~e.~$M:=\Hom_\mathbb{Z}(N,\mathbb{Z})$ and $\left<\cdot,\cdot\right>: M\times N\to \mathbb{Z}$ the canonical pairing. We define an algebraic torus $T_N$ by setting
$$T_N:=N\otimes_\mathbb{Z}\mathbb{C}^*=\Hom_{\mathbb{Z}}(M,\mathbb{C}^*),$$ then $T_N$ is isomorphic to $\mathbb{C}^*$ and $N$ can be considered as its group of one-parameter subgroups, while $M$ is its group of characters. An element $n\in N$ gives rise to a one-parameter subgroup $\gamma_n:\mathbb{C}^*\to T_N$ by setting $$\gamma_n(\lambda)(m):=\lambda^{\left<m,n\right>} \text{ for all } \lambda \in \mathbb{C}^*\text{ and }m\in M$$
and an element $m\in M$ defines a charakter $\chi(m):T_N\to \mathbb{C}^*$ by $$\chi(m)(t):=t(m)\text{ for all } t\in T_N.$$
Let $N_\mathbb{R}$ and $M_\mathbb{R}$ be the mutually dual $\mathbb{R}$-vector spaces $$N_\mathbb{R}:=N\otimes_\mathbb{Z}\mathbb{R}\text{ and }M_\mathbb{R}:=M\otimes_\mathbb{Z}\mathbb{R}=\Hom_\mathbb{Z}(N,\mathbb{R}).$$ 
\begin{defin} A convex polyhedral cone $\sigma$ is a subset $\sigma\subset N_\mathbb{R}$ such that there exist elements $n_1,\ldots,n_s\in N_\mathbb{R}$ with $$\sigma=\mathbb{R}_{\geq 0}n_1+\cdots + \mathbb{R}_{\geq 0}n_s.$$ The cone $\sigma$ is called rational, if $n_1,\ldots,n_s$ can be chosen to lie in $N$. The dual cone $\sigma^\vee\subset M_\mathbb{R}$ is defined by $$\sigma^\vee:=\{m\in M_\mathbb{R}\mid \left<m,n\right>\geq 0,\text{for all } n\in \sigma\}.$$ A subset $\tau\subset \sigma$ is called a face of $\sigma$, if there exists an element $m_0\in M_\mathbb{R}$ such that $m_0\in \sigma^\vee$ and $$\tau=\sigma\cap\{m_0\}^\perp:=\{n\in\sigma\mid \left<m_0,n\right>=0\}.$$
A fan $\Delta$ is a family of convex polyhedral cones in $N_\mathbb{R}$ subject to following conditions:
\renewcommand{\labelenumi}{(\roman{enumi})}
\begin{enumerate}\item Each $\sigma \in \Delta$ is rational.
\item Each $\sigma\in\Delta$ is strongly convex, meaning that $\sigma\cap(-\sigma)=\{0\}$.
\item Every face of $\sigma$ belongs to $\Delta$.
\item For $\sigma,\sigma'\in\Delta$, the intersection $\sigma\cap\sigma'$ is a face of $\sigma$ and $\sigma'$.
\end{enumerate}
\end{defin}
We now come to the fundamental result in this theory:
\begin{thrm}Via a covariant functor $(N,\Delta) \to T_N\emb(\Delta)$ the category of fans is equivalent to the category of toric varieties over $\mathbb{C}$.
\end{thrm}
\begin{remark}Here $T_N\emb(\Delta)$ is defined as the union of certain sets $U_\sigma$:$$T_N\emb(\Delta):=\bigcup_{\sigma\in \Delta}U_{\sigma},\text{ where }$$
$$U_\sigma:=\{u:M\cap\sigma^\vee\to \mathbb{C}\mid u(0)=1, u(m+m')=u(m)u(m')\text{ for } m,m'\in M\cap\sigma^\vee\}.$$ One can prove that these sets can be identified with algebraic subsets of $\mathbb{C}^n$ and the conditions (i)-(iv) in the definition of a fan guarantee that these sets glue together to define a toric variety.
\end{remark}



\subsection{The symplectic geometry of Hermitian symmetric spaces}
Let $(M,\omega)$ be a symplectic manifold and $H$ a compact connected Lie group acting on $M$ by symplectomorphisms, i.~e.~$g^*\omega=\omega$ for all $g\in H$. The action is called Hamiltonian if there exists a momentum map $\mu:M \to \frak h^*$, which is by definition an $H$-equivariant map such that $$\iota_{X_{\xi}}\omega=d\mu^\xi \text{ for all }\xi\in\frak h,$$ where $\mu^\xi:M\to \mathbb{R}$ is defined by $\mu^\xi(x):=\mu(x)(\xi)$. Here $X_{\xi}$ denotes the vector field on $M$ induced by $\xi \in \frak h$ and $\iota_{X_{\xi}}\omega$ the contraction of $\omega$ by $X_{\xi}$. Let $x\in M$ and $d\mu_x: T_xM\to T_{\mu(x)}\frak h^*$ be the differential of $\mu$. We then have $$\ker d\mu_x=(T_xH\kdot x)^{\perp_\omega}:=\{v\in T_xM\mid \omega(v,w)=0 \text{ for all }w\in T_xH\kdot x\}.$$
In general the existence of a momentum map is not clear, however there are some special situations, where the existence of $\mu$ can be guaranteed. 

Let $M$ be a coadjoint orbit for a  Lie group $G$. Then the momentum map for $G$ is just the inclusion of $M$ into $\frak g^*$. Any action of a Lie subgroup $H\subset G$ is again Hamiltonian and the momentum map for the $H$-action is just the projection on $\frak h^*$. Another important situation is the following: Let $H$ be a compact Lie group acting in a Hamiltonian way on $M$ and $N\subset M$ a submanifold which is invariant under $H$. Then the $H$-action on $M$ is again Hamiltonian and the corresponding momentum map is just the restriction of the original momentum map to $N$.

\noindent Let $N\subset M$ be a submanifold. We call $N$ a \emph{coisotropic} submanifold if \begin{equation}\label{a}T_pN^{\perp_\omega}\subset T_pN\end{equation} for all points $p\in N$. Let $J\in \End(TM)$ be a compatible almost complex structure, i.~e.~$\omega(JX,JY)=\omega(X,Y)$  and $g(X,Y):=\omega(X,JY)$ defines a Riemannian metric on $M$. Then condition (\ref{a}) can be equivalently formulated as \begin{equation}\nu_pN \perp_{g} J\nu_pN\end{equation} for all $p\in N$, where $\nu_pN$ denotes the normal space of $N$ in $M$, i.~e.~the normal bundle of $N$ is \emph{totally real}. We call a Lie group action on $M$ coisotropic, if all principal orbits are coisotropic.

Let $M$ be now a Hermitian symmetric space of compact type. Then any  isometric action on $M$ is Hamiltonian, since $M$ can be realized as an adjoint orbit.
The following theorem links the theory of polar actions to that of coisotropic ones.
\begin{thrm}[PoTh 2002] \label{p=c}A polar action on an irreducible compact homogeneous simply-connected K\"ahler manifold is coisotropic.
\end{thrm}
\noindent So in other words the sections are totally real submanifolds of $M$.

Another important theorem that characterizes coisotropic actions on compact K\"ahler manifolds is the following Equivalence Theorem of  Huckleberry and\break  Wurzbacher \cite{HuWu}.
\begin{thrm}[HuWu 1990]\label{equiv} Let $M$ be a connected compact K\"ahler manifold with an isometric and Hamiltonian action of a connected compact group $H$. Then the following conditions are equivalent:
\begin{enumerate} \renewcommand{\labelenumi}{(\roman{enumi})}
\item The $H$-action is coisotropic. \item The momentum map $\mu: M\to \frak h^*$ separates orbits. \item The cohomogoneity of the $H$-action is equal to the difference between the rank of $H$ and the rank of a regular isotropy subgroup of $H$. 
\item The space $C^{\infty}(M)^H$ of $H$-invariant functions on $M$ is abelian with respect to the Poisson bracket. \item The K\"ahler manifold $M$ is projective algebraic, $H^{\mathbb{C}}$-almost homogeneous and a spherical embedding of the open $H^{\mathbb{C}}$-orbit.
\end{enumerate}
\end{thrm}
\begin{remark} The points (ii) and (iii) will be of special importance for our work.
\end{remark}
\subsection{Proof of the Reduction Theorem}

For a polar action on a Riemannian manifold the normal distribution to the principal orbits is integrable, a leaf is given by the intersection of a section with the set of regular points of the action.  Now on a K\"ahler manifold the main result of \cite{Podesta} shows that the sections are totally real submanifolds, so it makes sense to define the following notion:
\begin{defin} Let $M$ be an irreducible compact homogeneous simply-con\-nected K\"ahler manifold and $H$ a compact connected Lie group acting polarly on $M$. We call the $H$-action complex integrable if and only if the distribution $$\mathcal{D}_q:=T_q\Sigma_q\oplus JT_q\Sigma_q, \;\; q\in M_r$$ is integrable. Here $M_r$ denotes the set of regular points of the $H$-action and $\Sigma_q$ is the unique section through $q$. 
\end{defin}
The first observation concerning a complex integrable polar action is the following:
\begin{lemma}\label{geodesic} Let $H\times M\to M$ be a complex integrable polar action. Then the leaves of the distribution $\mathcal{D}$ are automatically totally geodesic.
\end{lemma}
\begin{proof}Let $X,Y$ be two local sections of $\mathcal{D}$. We have to prove that the normal part of $\nabla_XY$ with respect to $\mathcal{D}$ vanishes. Here $\nabla$ denotes the Levi-Civita connection of $M$. We write $X=X_1+JX_2$ and $Y=Y_1+JY_2$ with $X_i,Y_i$ sections of $\mathcal{D}$ such that $X_i(q),Y_i(q)\in T_q\Sigma_q$. We then get
\begin{align*}\nabla_XY&=\nabla_{X_1}Y_1+\nabla_{X_1}JY_2+\nabla_{JX_2}Y_1+\nabla_{JX_2}JY_2\\&=\nabla_{X_1}Y_1+J\nabla_{X_1}Y_2+[JX_2,Y_1]+J\nabla_{Y_1}X_2+J([JX_2,Y_2]+J\nabla_{Y_2}X_2) \end{align*} which shows that $\nabla_XY$ is tangent to $\mathcal{D}$ since the sections are totally geodesic and the distribution $\mathcal{D}$ is assumed to be integrable.
\end{proof}
Examples of complex integrable polar actions are provided by polar actions which are complex asystatic. Recall that an isometric Lie group action is called asystatic if the isotropy representation of a principal isotropy subgroup has no non-trivial fixed vectors on the tangent space of the corresponding orbit. In \cite{Aleks3} it is proved that asystatic actions are polar and a section is given by the fixed point set of a principal isotropy subgroup.
Now on a K\"ahler manifold there is the related concept of a $\mathbb{C}$-asystatic action. This concept is introduced first in \cite{Podesta} and systematically studied in  \cite{Gori2}, where such actions are called symplectically asystatic, but for polar actions we prefer the notion $\mathbb{C}$-asystatic as we will explain.  The following result is proved in \cite{Gori2}:
\begin{prop}Let $M$ be K\"ahler manifold endowed with an effective Hamiltonian action of a compact group $K$. Let $(L)$ be a principal isotropy type, $\overline{M}$ a core and $c$ the cohomogeneity of the $K$-action. Then \begin{enumerate}\renewcommand{\labelenumi}{(\roman{enumi})}\item $\dim\overline{M}\geq c+{\rk}\; K-\rk\; L$;
\item the equality in (i) holds if and only if the group $(N_K(L)/L)^\circ$ is abelian.
\end{enumerate} 

\end{prop}
The \emph{core} of the $K$-action is defined as the closure of $M^L\cap M_r$, where $M^L$ denotes the fixed point set of $L$.
Now Gori and Podest\`a call the $K$-action \emph{symplectically asystatic} if in $(i)$ equality holds. 
\vskip 10pt
Let $H$ be now a compact connected Lie group acting polarly on an irreducible compact homogeneous simply-connected K\"ahler manifold. Note that for this action, which is coisotropic, $(i)$ reduces to $$\dim\overline{M}\geq 2c.$$
Let $p\in M$ be a regular point and $\Sigma$ the section through $p$. Then the tangent space of the $H$-orbit through $p$ splits as $$T_pH\kdot p=JT_p\Sigma\oplus (T_pH\kdot p\cap JT_pH\kdot p).$$
Assume now that in (i) equality holds. Then this means that the isotropy $H_p$ has no non-trivial fixed vector on the maximal holomorphic subspace $T_pH\kdot p\cap JT_pH\kdot p$ of $T_pH\kdot p$, so we define:
\begin{defin} The $H$-action is called $\mathbb{C}$-asystatic if a principal isotropy subgroup has no non-trivial fixed vector on the maximal holomorphic subspace of the tangent space of the corresponding orbit.
\end{defin} 

The following is an example of a polar $\mathbb{C}$-asystatic action (cf.~\cite{Gori}): Let $M=G/K$ be the complex quadric $Q_n=\SO(n+2)/\SO(2)\times \SO(n)$. Then the $K$-action on $M$ is $\mathbb{C}$-asystatic. A principal isotropy subgroup $L$ is given by \break $\mathbb{Z}_2\times \SO(n-2)$ and $(N_K(L)/L)^\circ=\SO(2)^2$ is abelian.
\vskip 10pt
\noindent We are interested in $\mathbb{C}$-asystatic actions because of the following proposition:
\begin{prop} A polar $\mathbb{C}$-asystatic action is complex integrable and  a leaf $S_p$ through a point $p \in M_r$ is given by the intersection of $M_r$ with the connected component of the fixed point set of $H_p$ which contains $p$.
\end{prop}

\begin{proof} Let $S_p$ denote the intersection of $M_r$ with the connected component of the fixed point set of $H_p$ which contains $p$. Since $p$ is regular, $H_p$ acts trivially on $T_p\Sigma_p \oplus JT_p\Sigma$ and hence $T_pS_p\supset T_p\Sigma_p\oplus JT_p\Sigma_p$. On the other hand $T_pM$ splits as 
$$T_pM=T_pH\kdot p\oplus T_p\Sigma_p=JT_p\Sigma_p\oplus (T_pH\kdot p\cap JT_pH\kdot p) \oplus T_p\Sigma_p$$
and we have seen that $H_p$ has no fixed points on $T_pH\kdot p\cap JT_pH\kdot p$, so we actually have $T_pS_p=T_p\Sigma_p\oplus JT_p\Sigma_p$.\\
The next step is to show that we have $T_xS_p=T_x\Sigma_x\oplus JT_x\Sigma_x$ for an arbitrary point $x\in S_p$. Since $x$ is a regular point in the fixed point set of $H_p$, we have $H_x=H_p$. Using the same notations and arguments as above we have $T_xS_x=T_x\Sigma_x\oplus JT_x\Sigma_x$, but $M^{H_x}=M^{H_p}$ and hence $S_x=S_p$ and the claim follows.
\end{proof}

From now on we assume $M=G/K$ to be an irreducible Hermitian symmetric space of compact type and $H$ a connected compact Lie group. Moreover the $H$-action on $M$ is assumed to be polar and complex integrable. Let $\mu: M \to \frak h$ be a momentum map for the $H$-action, which clearly exists, since $M$ can be realized as the coadjoint orbit of the full isometry group of $M$. Since $H$ is compact, we have identified $\frak h$ with its dual.
\noindent We start with an elementary lemma:
\begin{lemma} Let $p \in M$ be a regular point for the $H$-action and $\Sigma$ the section through $p$. Then $\ker(d\mu_p)=JT_p\Sigma$. 
\end{lemma}
\begin{proof} First of all since $p$ is regular $T_p\Sigma$ is precisely the normal space of the orbit $H\kdot p$. The lemma follows then  immediately from the standard fact that $\ker(d\mu_p)=(T_pH\kdot p)^{\perp_\omega}$, where $ {\perp_\omega}$ denotes the orthogonal complement with respect to the K\"ahler form $\omega$ of $M$. 
\end{proof}

\noindent Now define $$M_{\mu}:=\{x\in M \mid \dim H \kdot\mu(x)\geq  \dim H\kdot\mu(y)\text{ for all }y\in M\},$$ then in \cite{HuWu} it is proved that the set $M_{\mu}$ is open and that for all $x \in M_{\mu} $ the restricted momentum map
$$\mu :H\kdot x \longrightarrow  H\kdot \mu(x)$$ is a  principal $(S^1)^{r(x)}$-bundle, where $r(x):=\dim\,(T_xH\kdot x)\cap (T_xH\kdot x)^{\perp_\omega}$. For a polar action we have the following refinement:
\begin{lemma}\label{torus} For any regular point $p\in M$ of the $H$-action $T:=H_{\mu(p)}/H_p$ is a torus with $\dim T=\dim \Sigma$.
\end{lemma}
\begin{proof} Since the $H$-action is coisotropic we have $T_pH\kdot p^{\perp_\omega}\subset T_pH\kdot p$ and hence $r(p)=\dim\, \Sigma$ for all regular points $p$. Hence it is enough to show that every regular point of the $H$-action is contained in $M_\mu$. Since the set $M_\mu$ is open, it is enough to show that $\dim\,H\kdot \mu(p)$ is constant on the regular points of the $H$-action. For a regular point $p \in M$ we have $$T_pH\kdot p=JT_p\Sigma \oplus (T_pH\kdot p\cap JT_pH\kdot p)$$ and $T_pH\kdot p\cap JT_pH\kdot p$ is isomorphic to $T_{\mu(p)}H\kdot \mu(p)$ via $d\mu_p$, which proves the lemma.
\end{proof}
Now the torus $T=H_{\mu(p)}/H_p$ is at first an abstract torus and it is not clear how it acts on $M$. However on the Lie algebra level we have
$$\frak h_{\mu(p)}=\frak t\oplus \frak h_p,$$ with $\frak t$ abelian. So in the following we mean by defintion that the notation $T=H_{\mu(p)}/H_p$ is a symbol which stands for the torus which is defined by $T=\exp\;\frak t$. So $T$ is a subgroup of $H_{\mu(p)}$ and hence clearly acts on $M$.

\begin{lemma} Let $p \in M$ be a regular point and  $T:=H_{\mu(p)}/H_p$. Then $T\kdot p$ describes the $\mu$-fiber through $p$, i.~e.~$\mu^{-1}(\mu(p))=T\kdot p$. Moreover we have $T_pT\kdot p=JT_p\Sigma$.
\end{lemma}
\begin{proof}Viewing the restricted momentum map $$\mu: H\kdot p \to H\kdot \mu(p)$$ as the homogeneous fibration $$H/H_p \to H/H_{\mu(p)}$$ we see that $H_{\mu(p)}\kdot p$ is the fiber of the restricted momentum map. On the other hand we know that $\mu$ separates orbits since the action is coisotropic and hence $H_{\mu(p)}\kdot p$ is the fiber of the full momentum map. By the definition of $T$ it follows immediately that $T\kdot p=H_{\mu(p)}\kdot p$: Let $g\in H_{\mu(p)}$. Since $H_{\mu(p)}$ is the centralizer of a torus it is connected and hence there exists $X\in \frak h_{\mu(p)}$ with $\exp(X)=g$. Write $X=X_1+X_2$ with $X_1\in \frak t$ and $X_2\in \frak h_p$. Since $X_1$ and $X_2$ commute we have $$gp=\exp(X)p=\exp(X_1)\exp(X_2)p=\exp(X_1)p\in T\kdot p.$$The last claim is clear from the equality  $\ker(d\mu_p)=JT_p\Sigma$. 
\end{proof}

Recall the distribution $\mathcal{D}$ defined by $$ \mathcal{D}_q:=T_q\Sigma_q\oplus JT_q\Sigma_q,$$
where $q\in M_r$ and $\Sigma_q$ is the unique section through $q.$ This gives us a well defined distribution on $M_r$ and by assumption this distribution is integrable. We have seen in Lemma \ref{geodesic} that the leaves are totally geodesic although they are not complete, since the distribution $\mathcal{D}$ is only defined on the regular set of the $H$-action. However since the ambient space is a symmetric space we can easily define a completion of a leaf $S_p$. Let $\Sigma$ be the section through $p$. Without loss of generality we can assume that $p$ is the origin of $M=G/K$. Let $\frak g=\frak k\oplus\frak p$ be the corresponding Cartan decomposition and $\frak m\subset \frak p$ the tangent space of $\Sigma$. Then Lemma \ref{geodesic} implies that $\frak m \oplus J\frak m$ is a Lie triple system. Hence the following definition makes sense:

\begin{defin} The complexification $\Sigma^{\mathbb{C}}\subset M$ of a section $\Sigma$ is defined by $$\Sigma^{\mathbb{C}}:=\exp_p(T_p\Sigma\oplus JT_p\Sigma),$$ where $p\in \Sigma$ is any regular point. \end{defin}
The following proposition is then clear:
\begin{prop} The complexification $\Sigma^{\mathbb{C}}$ is a totally geodesic complex submanifold of $M$, which contains $\Sigma$. Moreover if $p\in \Sigma^{\mathbb{C}}$ is a regular point, then the leaf of $\mathcal{D}$ through $p$ is an open subset of $\Sigma^{\mathbb{C}}$. If $q\in \Sigma^{\mathbb{C}}$ is any regular point and $\Sigma_q$ is the section through $q$ then $\Sigma^{\mathbb{C}}$ is also the complexification of $\Sigma_q$ and we have $\Sigma^{\mathbb{C}}=\exp_q(T_q\Sigma_q\oplus JT_q\Sigma_q)$. 
\end{prop}
\begin{remark} The next proposition will show that alternatively one can define $\Sigma^{\mathbb{C}}:=\overline{T^{\mathbb{C}}\kdot p}$, where $p\in \Sigma$ is a regular point and $T=H_{\mu(p)}/H_p$. This is interesting in its own, since in general the closure of complex torus orbits can be  highly singular.
\end{remark}
\begin{prop} Let $p\in \Sigma$ be a regular point of the $H$-action on $M$ and $T=H_{\mu(p)}/H_p$. Then $T$ acts effectively, isometrically and in a Hamiltonian way on $\Sigma^{\mathbb{C}}$.
\end{prop}
\begin{proof} We will show that $H_{\mu(p)}$ acts on $\Sigma^{\mathbb{C}}$ and that $H_p$ is the kernel of this action. Remembering the homogeneous fibration $$H/H_{p}\to H/H_{\mu(p)}$$ we see that $T_pH_{\mu(p)}\kdot p=JT_p\Sigma_p$ and hence it follows for any $g\in H_{\mu(p)}$ that 
$$T_{gp}H_{\mu(p)}\kdot p=dg(T_pH_{\mu(p)}\kdot p)=dg(JT_p\Sigma_p)=JT_{gp}\Sigma_{gp}\subset \mathcal{D}_{gp},$$  which shows that $H_{\mu(p)}\kdot p\subset S_p\subset \Sigma^{\mathbb{C}}$. Now we can show that for an arbitrary point $x\in \Sigma^{\mathbb{C}}$ the orbit $H_{\mu(p)}\kdot x$ lies again in $\Sigma^{\mathbb{C}}$. Write $x=\exp_p(v+Jw)$ with $v,w \in T_p\Sigma$. Using the fact that $H_{\mu(p)}$ acts isometrically we get for $g\in H_{\mu(p)}$: $$gx=g\,\exp_p(v+Jw)=\exp_{gp}(dgv+Jdgw)\in \Sigma^\mathbb{C},$$
since we already know that $gp\in \Sigma^{\mathbb{C}}$ and $ dgv+Jdgw \in  T_{gp}\Sigma^\mathbb{C}$.

\noindent It is easy to see that the $T$-action on $\Sigma^{\mathbb{C}}$ is effective. We have to show that $H_{p}=\bigcap_{q\in\Sigma^\mathbb{C}}(H_{\mu(p)})_q$. The inclusion "$\supset$" follows from the trivial observation $(H_{\mu(p)})_p=H_p$ while "$\subset$" follows from the equality $\Sigma^\mathbb{C}=\exp_p(T_p\Sigma\oplus JT_p\Sigma)$ and the triviality of the slice representation  of $H_p$ on $T_p\Sigma$. It is only left to prove that the action is Hamiltonian. But this is clear since the $T$-action on $M$ is Hamiltonian and $T$ leaves $\Sigma^\mathbb{C}$ invariant.
\end{proof}
\noindent For the final proof we will need the following criterion for polarity:
\begin{prop}[\cite{Go}] Let $M=G/K$ be a symmetric space of compact type with
a Riemannian metric induced from some $\Ad(G)$-invariant inner product on the Lie
algebra $\frak g$ of G. Consider a closed, connected subgroup $H\subset G$. By replacing $H$ by a
conjugate, if necessary, we may assume that  $eK$ is a regular point. Write
$\frak g=\frak k \oplus \frak p$ for the Cartan decomposition, denote by $\frak h$ the Lie algebra of $H$, and define $\frak m=\frak p\cap \frak h^\perp$ $($i.~e.~$\frak m$ is the normal space to the $H$-orbit through $eK)$. Then the action of $H$ on $M$ is polar if and only if the following two conditions hold:\renewcommand{\labelenumi}{(\roman{enumi})}
\begin{enumerate}\item $[[\frak m,\frak m],\frak m]\subset \frak m$
\item $[\frak m, \frak m]\perp \frak h.$
\end{enumerate}
\end{prop}

\begin{proof} [Proof of the Reduction Theorem] Since $\frak h\supset \frak t$ the polarity of the $T$-action is clear using the criterion of Gorodski.
\noindent Since $T$ acts effectively and isometrically we have $T\subset \iso(\Sigma^\mathbb{C})$ and hence ${\rk}(\iso(\Sigma^\mathbb{C}))\geq \dim\,T=\dim_{\mathbb{C}}\Sigma^\mathbb{C}$. On the other hand we have the general fact that for an irreducible Hermitian symmetric space $N$ always $\dim_{\mathbb{C}}\,N\geq{\rk}(\iso(N))$ and equality holds if and only if $N$ is a complex projective space. This shows $\Sigma^\mathbb{C}=\mathbb{P}_{n_1}(\mathbb{C})\times \cdots \times \mathbb{P}_{n_k}(\mathbb{C})$. This argument was also used in the proof of Corollary 1.4. in \cite{Podesta3}.
\end{proof}
\subsection{Applications of the Reduction Theorem}
\subsubsection{A hyperpolarity criterion}
Let us return to the momentum map $\mu:M\to \frak h$ of the $H$-action. If $p\in M$ is a regular point of the $H$-action, then we have seen that the orbit $T\kdot p$ describes the momentum fiber through $p$. Interestingly these gives a criterion to distinguish a polar action from a hyperpolar one. More precisely:
\begin{prop}\label{criterion} Let $M$ be an irreducible Hermitian symmetric space endowed with a complex integrable polar action of a compact connected Lie group $H$. The $H$-action is hyperpolar if and only if there exists a totally geodesic momentum fiber through a regular point of the $H$-action. 
\end{prop}
First we will prove a lemma which describes a splitting of the induced torus action on the complexified section.
\begin{lemma}\label{split} Let $p\in\Sigma^\mathbb{C}=\mathbb{P}_{n_1}(\mathbb{C})\times \cdots \times \mathbb{P}_{n_k}(\mathbb{C})$ be a regular point of the $H$-action and $T:=H_{\mu(p)}/H_p$. Then $T$ splits into $T=(S^1)^{n_1}\times \cdots \times (S^1)^{n_k}$ such that the $T$-action is the product of polar actions of $(S^1)^{n_i}$ on $\mathbb{P}_{n_i}(\mathbb{C})$.
\end{lemma}
\begin{proof}This is trivial. One only has to observe that $T$ is a maximal torus in the isometry group of $\Sigma^\mathbb{C}$.
\end{proof}

The next lemma is also clear:
\begin{lemma} Let $\Sigma'$ be a section for the $H$-action which is contained in $\Sigma^\mathbb{C}$. Then $\Sigma'$ is also a section for the $T$-action on $\Sigma^\mathbb{C}$, where $T=H_{\mu(p)}/H_p$ for any regular point $p\in \Sigma^\mathbb{C}$. Conversely any section for the $T$-action on $\Sigma^\mathbb{C}$ is also a section for the $H$-action on $M$.
\end{lemma}
\begin{proof} Pick a regular point $x \in \Sigma'$. Then there exists a section $\hat{\Sigma}\subset \Sigma^\mathbb{C}$ through $x$ for the $T$-action. This means there exists an element $t\in T$ with $t\kdot\Sigma=\hat{\Sigma}$. Since $x$ is a regular point  and $T\subset H$ it follows that $\Sigma'=\hat{\Sigma}$. The same argument proves the other direction.
\end{proof}
In \cite{HuWu} it was shown that for all $p\in M_{\mu}$ the torus orbit $T\kdot p$ with $T:=H_{\mu(p)}/H_p$ describes the $\mu$-fiber through $p$. In our situation one can show that this torus describes not only the $\mu$-fiber through $p$, but also the fibers on a larger set. More precisely:
\begin{prop} 
Let $p\in\Sigma^\mathbb{C}\subset M$ be a regular point of the $H$-action and $T:=H_{\mu(p)}/H_p.$ Then $T\kdot x$ is the $\mu$-fiber through $x$ for all regular $x\in \Sigma^\mathbb{C}$.
\end{prop}
\begin{proof}Let $x\in \Sigma^\mathbb{C}\subset M$ be a regular point of the $H$-action and $T':=H_{\mu(x)}/H_x$. If $\Sigma_x$  is the section for the $H$-action through $x$ then it follows from the lemma above that $\Sigma_x$ is contained in $\Sigma^\mathbb{C}$ and that it is a section for the $T$-action as well as for the $T'$-action. It follows that the sections for the $T$- and the $T'$-action are identical and hence the principal orbits of both actions coincide.
\end{proof}
\begin{proof}[Proof of Proposition \ref{criterion}] Let us first assume that there exists a totally geo\-desic $\mu$-fiber through a regular point $p\in M$ of the $H$-action. Let $\Sigma$ be the section through $p$ and $T:=H_{\mu(p)}/H_p$. Then we know that $T$ acts polarly on $\Sigma^\mathbb{C}$ with $\Sigma$ as a section. Let $N:=T\kdot p\subset \Sigma^\mathbb{C}\subset M$ denote the totally geodesic $\mu$-fiber through $p$. Pick two vectors $u,v\in T_pN$ and let $R^{L}$ denote the curvature tensor of a submanifold $L\subset M$. If $g, J$ denote the metric and the complex structure of $M$ we have
\begin{align*}g(R^{\Sigma}(Ju,Jv)Jv,Ju)&=g(R^{\Sigma^\mathbb{C}}(Ju,Jv)Jv,Ju)\\&=g(R^{\Sigma^\mathbb{C}}(u,v)v,u)=g(R^N(u,v)v,u)=0,\end{align*}
since $N$ is a totally geodesic orbit of a torus and thus flat.

The converse direction follows from Lemma \ref{split}, since it shows that for a hyperpolar action the $T$-action splits into a product of $S^1$-action on a product of Riemann spheres.
\end{proof}
\subsubsection{The rank one case}
Let $X$ be a complex projective space. Then any totally geodesic submanifold of $X$ is again a rank one space, hence it is known, that a section for a polar action on $X$ must be again a rank one space. Thorbergsson and Podest\`a have proven in \cite{Podesta2} by classifying all polar actions on $X$ that  a section is isometric to a real projective space. This will also follow from the following proposition:
\begin{prop}Let $M$ be an irreducible Hermitian symmetric space endowed with a complex integrable polar action of a compact connected Lie group $H$. Let $\Sigma$ be a section for the $H$-action. Then $\Sigma$ is isometric to a product of real projective spaces, i.~e.~$\Sigma=\mathbb{P}_{n_1}(\mathbb{R})\times \cdots \times \mathbb{P}_{n_k}(\mathbb{R})$. The $H$-action is hyperpolar if and only if $n_1=\cdots=n_k=1$.
\end{prop}

\begin{proof} We have seen that for a section $\Sigma$ the complexification $\Sigma^\mathbb{C}$ is a product of complex projective spaces, i.~e.~$\Sigma^\mathbb{C}=\mathbb{P}_{n_1}(\mathbb{C})\times \cdots \times \mathbb{P}_{n_k}(\mathbb{C})$. Let $p\in \Sigma$ be a regular point of the $H$-action and $T=H_{\mu(p)}/H_p$.
We have seen in Lemma \ref{split} that the $T$-action splits into $(S^1)^{n_i}$-actions on $\mathbb{P}_{n_i}(\mathbb{C})$. These actions are clearly conjugated to the action of the standard maximal torus in the isometry group of $\mathbb{P}_{n_i}(\mathbb{C})$, for which it is well known that it is polar with a section being isometric to $\mathbb{P}_{n_i}(\mathbb{R})$.
\end{proof}
\begin{prop} Let $X$ be a complex projective space and $H$ a connected compact Lie group acting polarly on $X$. Then the action is complex integrable and a section is isometric to a real projective space.
\end{prop}
\begin{proof}We have to show that the distribution $\mathcal{D}$ is integrable. This follows from the fact that any complex subspace in the tangent space of $X$ is a Lie triple system (cf.~\cite{KN2} p.~277). It is then clear that we get a leaf by defining  $S_p:=\exp_p(T_p\Sigma_p\oplus JT_p\Sigma_p)\cap M_r$. The claim follows then from the proposition above since $X$ is a rank one space.
\end{proof}
\subsection{Hermitian symmetric spaces of maximal rank}

In this section we will analyze polar actions on the Lagrangian Grassmannian of maximal isotropic planes in $\mathbb{C}^{2n}$, i.~e.~the Hermitian symmetric space $M=G/K=\Sp(n)/\U(n)$. The special property of this space that makes it interesting is the fact that it is a space of maximal rank, i.~e.~the rank of $M$ is equal to the rank of its isometry group. This allows us to use some general results about toric varieties and torus actions on generalized flag manifolds. Recall the non-degenerate stratum $Z(\T)$ which was defined as the intersection of the big cells of all Borel subgroups containing a fixed maximal Torus $\T$. 
\begin{prop}\label{Z} Let $H\times M\to M$ be a complex integrable hyperpolar action with maximal cohomogeneity $($i.~e.~$\cohom(H,M)=\rk\; G)$ of a compact connected Lie group on $M$. Then every regular point $p\in M$ is contained in the non-degenerate stratum $Z(\T)$ of the complexification of the torus $T=H_{\mu(p)}/H_p$.
\end{prop}


The first observation concerning this stratum was made by Gel'fand and\break  Serganova in \cite{GeSe}. First they define the notion of a \emph{thin cell}. Let $\G$ be a complex semisimple Lie group and $P\subset \G$ a parabolic subgroup. Then $G^\mathbb{C}/P$ is a generalized flag manifold. Let $B\subset G^\mathbb{C}$ be a  Borel subgroup containing a fixed Cartan subgroup $\T$. The orbits of $B$ in $M$ are called the Schubert cells associated with $B$.
\begin{defin} For every Borel subgroup $B$ containing $\T$ we choose precisely one Schubert cell. The intersection of all these chosen Schubert cells is called a thin cell, if it is not empty.
\end{defin}
\begin{remark} Let $W=N_G(\T)/\T$ denote the Weyl group. Then $W$ acts simply transitively on the set of all Borel subgroups containig $\T$. Let $B$ be a fixed Borel subgroup containing $\T$ and $w\in W$. Choosing for the Borel subgroup $wBw^{-1}$ the Schubert cell through $wP$ we see that $Z(\T)$ is a thin cell.
\end{remark}
Let $T$ be a compact real form of $\T$ and $\mu_T: M \to \frak t$ the momentum map for the action of $T$ on $M$. It is well known that the image of $\overline{\T\kdot p}$ under the momentum map is a convex polytope for all $p\in M$ \cite{Atiyah}, \cite{Guillemin}. We have the following result of Gel'fand and Serganova:
\begin{thrm}[Theorem 1 in \S\  6 of \cite{GeSe}]\label{Gel} For two points $p,q\in M$ we have $\mu_T(\overline{\T\kdot p})=\mu_T(\overline{\T\kdot q})$ if and only if $p$ and $q$ lie in the same thin cell.
\end{thrm}
\begin{remark} It is important to note that this theorem does not imply that $\overline{\T\kdot p}$ and $\overline{\T\kdot q}$ are isomorphic as toric varieties. The general theory of toric varieties is not directly applicable, because nothing is said about the normality of these varieties. This is a mistake Flaschka and Haine make in {\rm \cite{FlHa}}. Their results only hold for those torus orbit closures in $G^\mathbb{C}/P$ which are normal. In particular contrary to what is claimed by them, there exist torus orbit closures in $Z(\T)$ which are not isomorphic.
\end{remark}

However one can show that there exists an open subset of the non-degenerate stratum $Z(\T)$ such that all torus orbit closures are normal and isomorphic as toric varieties, as was done by Dabrowski \cite{Da}. Let us call this set $Y$ and torus orbit closures of elements in $Y$ generic. We will not go into the detailed description of this set, since we do not need it. We just briefly state the main result of \cite{Da}:
\begin{thrm}[Theorem 3.2 of \cite{Da}]\label{Dab1}Let $x\in Y\subset G^\mathbb{C}/P$ and \break $X:=\overline{\T x}$. Then: \renewcommand{\labelenumi}{(\roman{enumi})}\begin{enumerate}\item $X$ is a normal variety.
\item The fan corresponding to $X$ consists of the cones 
$$C_w=-w\bigcup_{z\in W_P} zD,\;\; w\in W^P$$ together with their faces. In particular the closure of any two orbits in $Y$ are isomorphic as $\T$-equivariant embeddings.
\end{enumerate}
Here $D$ denotes the fundamental Weyl chamber with respect to a chosen basis of simple roots, $W_P$ is the subgroup of the Weyl group generated by the simple roots defining $P$ and $W^P$ is a fixed set of representatives of the coset space $W/W_P$.
\end{thrm}
The next step will be to identify the generic torus orbit type for $M=\Sp(n)/\U(n)$.  As we have seen, this can be done by identifying the fan describing this toric variety. 

Let $\Delta:=\{\pm(e_k\pm e_j)\mid 1\leq j,k\leq n, j\neq k\}\cup\{2e_i\mid 1\leq i\leq n\}$ be the root system of $\frak g^\mathbb{C}=\frak sp(n,\mathbb{C})$. Then a set of simple roots is given by $\pi=\{\alpha_1,\ldots,\alpha_n\}$ with  $\alpha_j=e_j-e_{j+1}$ for $j=1,\ldots,n-1$ and $\alpha_n=2e_n$. The parabolic subgroup $P\subset G^\mathbb{C}$ associated to $\pi_P:=\{\alpha_1,\ldots,\alpha_{n-1}\}$ realizes $M$ as the quotient $G^\mathbb{C}/P$. We realize this root system in an euclidean space $(E, (\cdot,\cdot))$ and denote by $W$  the Weyl group of it and by $W_P$ the subgroup of $W$ generated by reflections $s_\alpha$ with $s_\alpha \in \pi_P$. Let $D$ denote the fundamental Weyl chamber \break$\{v\in E\mid (v,\alpha)\geq 0 \text{ for all }\alpha \in \pi \}$. We have the following lemma:
\begin{lemma}\label{cone} $\bigcup_{z\in W_p}zD=\{v\in E \mid (v,e_i)\geq 0, i=1,\ldots,n\}$
\end{lemma}
\begin{proof}Let $L$ denote the right hand side. The first step will be to show that $D\subset L$. So let $v\in D$ then $(e_j-e_{j+1},v)\geq 0$ for all $j\in \{1,\ldots,n-1\}$ and $(e_n,v)\geq 0$. Setting $j=n-1$ and looking at the inequalities $(e_{n-1}-e_n,v)\geq 0$ and  $(e_n,v)\geq 0$ we see that $(e_{n-1},v)\geq 0$. Inductively we obtain $(e_i,v)\geq 0$ for all $i\in \{1,\ldots,n\}$, so $D\subset L$. Now we will show that the reflections $s_\alpha$ with $\alpha\in \pi_P$ leave $L$ invariant. Recall that the reflection $s_\alpha$ is given by $s_\alpha(w)=w-\frac{2(w,\alpha)}{(\alpha,\alpha)}\alpha$, we then obtain for $w\in L$ and $i\in \{1,\ldots,n-1\}, j\in \{1,\ldots,n\}$:
$$(s_{e_i-e_{i+1}}(w),e_j)=(w,s_{e_i-e_{i+1}}(e_j))=\left\{\begin{array}{ll} (w,e_{i+1}), &j=i \\ (w,e_i),&j=i+1\\(w,e_j),& j\notin \{i,i+1\}\end{array}\right\}\geq 0,$$ so $s_{e_i-e_{i+1}}(w)\in L$. Since the reflections $s_\alpha$ with $\alpha\in \pi_P$ generate $W_P$, we see that $W_P$ leaves $L$ invariant and hence $\bigcup_{z\in W_P}zD\subset L$.

Now let $v\in L$. Then $v=\sum_{i=1}^{n}a_ie_i$ with $a_i\in \mathbb{R}_{\geq 0}$ and the reflection $s_{\alpha_j}=s_{e_j-e_{j+1}}$ acts on $v$ by permuting the coefficients $a_j$ and $a_{j+1}$. It is then clear that there exist numbers $i_1,\ldots,i_k\in \{1,\ldots, n-1\}$ sucht that \break$w:=s_{\alpha_{i_1}}\cdots s_{\alpha_{i_k}}(v)=\sum_{i=1}^{n}c_ie_i$ with $c_1\geq c_2\geq\ldots\geq c_n\geq 0$. It follows $(w,e_i-e_{i+1})=c_i-c_{i+1}\geq 0$ for all $i\in \{1,\ldots,n-1\}$, i.~e.~$w\in D$ and hence $v\in zD$ with $z=(s_{\alpha_{i_1}}\cdots s_{\alpha_{i_k}})^{-1}$. Since $z\in W_P$ the proof of the lemma is complete.
\end{proof}
\begin{prop}\label{gen} The generic torus orbit type for $M$ is isomorphic to the product of $n$ complex projective lines.
\end{prop}
\begin{proof} Let $Q_{\varepsilon_1,\ldots,\varepsilon_n}:=\{\sum_{i=1}^{n}\varepsilon_ia_ie_i \mid a_i\in \mathbb{R}_{\geq 0}\}$ with $\varepsilon_i\in \{1,-1\}$ and $\mathcal{Q}$ the set of all such $Q_{\varepsilon_1,\ldots,\varepsilon_n}$. For the proof of this proposition it will be enough to show that \begin{equation}\label{2}\mathcal{Q}=\{-w\bigcup_{z\in W_P} zD\mid w\in W\}.\end{equation} First of all we note that thanks to Lemma \ref{cone} $\bigcup_{z\in W_p}zD=Q_{1,\ldots,1}$. As already seen the  reflection $s_{\alpha_j}$ with $1\leq j\leq n-1$ acts on $v=\sum_{i=1}^{n}a_ie_i$ by permuting the coefficients $a_j$ and $a_{j+1}$. For the reflection $s_{\alpha_n}$ we have $$s_{\alpha_n}(e_i)=\left\{\begin{array}{ll}\;\;\;e_i & i\neq n \\
-e_n &i=n \end{array}\right.$$ and hence $w\kdot Q_{1,\ldots,1}=Q_{\varepsilon_1,\ldots,\varepsilon_n}$ for any $w\in W$ and some $\varepsilon_i \in \{1,-1\}$. It follows that the right hand side of (\ref{2}) is contained in the left hand side. Conversely we can create any $Q_{\varepsilon_1,\ldots,\varepsilon_n}\in \mathcal{Q}$ starting with $Q_{1,\ldots,1}$ by successive applications of $s_{\alpha_n}$ (in order to create a minus sign) and the $s_{\alpha_j}$ with $1\leq j\leq n-1$  (in order to bring the minus sign at the right place). Now it is easy to see that $\mathcal{Q}$ describes indeed the fan of $\P_1(\mathbb{C})^n$.
\end{proof}
We will now prove Proposition \ref{Z}.
\begin{proof}[Proof of Proposition \ref{Z}] 

Let $p\in M$ be a regular point of the $H$-action. It follows from the Reduction Theorem that the complexification of the section through $p$ is isomorphic to a product of complex projective lines and hence $\overline{\T\kdot p}=\P_1(\mathbb{C})^n$. Let $\mu_T: M\to \frak t^*$ denote the momentum map for the $T$-action on $M$. Let $L^T$ denote the set of $T$-fixed points in $L$ for a subset $L\subset M$. Let $q\in Z(\T)\subset M$ be a point that realizes the generic torus orbit type of $M$. First we observe that the Euler characteristic of $M$ is $2^n$ (see for example \cite{Sanchez}) and hence both $(\overline{\T\kdot p})^T$ and  $(\overline{\T\kdot q})^T$ contain all $T$-fixed points in $M$. Let $\conv(R)$ denote the convex hull of a subset $R\subset \frak t$. Using the convexity results of \cite{Atiyah}, \cite{Guillemin}, \cite{Kirwan} we obtain
$$ \mu_T(\overline{\T\kdot p})=\conv(\mu ((\overline{\T\kdot p})^T) )=\conv(M^T)=\conv(\mu ((\overline{\T\kdot q})^T) )= \mu_T(\overline{\T\kdot q}),$$ i.~e.~the momentum images of the torus orbit closures through $p$ and $q$ are identical and hence it follows with Theorem \ref{Gel} that $p\in Z(\T)$ and the proposition is proven.
\end{proof}

\section{Polar actions on Wolf spaces}
Let us briefly recall the basic results about quaternion-K\"ahler manifolds, mainly to fix notations. We omit all proofs that can be found easily in the literature.

Let $\mathbb{H}^n$ denote the quaternionic $n$-space. We can identify $\mathbb{H}^n$ with $\mathbb{R}^{4n}$ such that the standard 
Euclidean inner product is given by $\left<u,v\right>={\rm Re}(v^*w)$, where ${}^*$ denotes transposing and quaternionically conjugating the entries. We can then define $\Sp(n):=\{A\in \GL(n,\mathbb{H})\mid \left<Au,Av\right> \text{ for all } u,v\in \mathbb{H}^n\}$. It is easily seen that the action of $\Sp(n)$ on $\mathbb{R}^{4n}$ is isometric, so $\Sp(n)\subset \SO(4n)$. However $\Sp(n)$ is not a maximal subgroup of $\SO(4n)$ since it commutes with the right multiplication with $\Sp(1)$. One then defines $\Sp(n)\kdot\Sp(1):=\Sp(n)\times\Sp(1)/\mathbb{Z}_2$, where $\mathbb{Z}_2$ is generated by $(-Id,-1)$ and this group is indeed a maximal subgroup of $\SO(4n)$ (\cite{Gray}). Now $\Sp(n)\kdot \Sp(1)$ appears also in the holonomy list of Berger so it is completely natural to define:
\begin{defin} A 4n-dimensional $(n>1)$ Riemannian manifold is called qua\- ternion-K\"ahler if its holonomy group is contained in $\Sp(n)\kdot \Sp(1)$. If $n=1$ we additionally require $M$ to be Einstein and self-dual.
\end{defin}
There exists the following equivalent definition from a slightly different viewpoint:
\begin{defin}\label{def2} Let $(M,g)$ be a Riemannian manifold. We call $M$ a quaternion-K\"ahler manifold if there exists a parallel rank three subbundle $Q\subset \End(TM)$ which is locally spanned by three almost complex structures $I,J,K$ satisfying the usual quaternionic relations.
\end{defin}
An important tool in the analysis of quaternion-K\"ahler manifolds is the so called $\emph{twistor space}$. Consider the unit-sphere subbundle $Z$ of the quaternionic bundle $Q\subset \End(TM)$. Let $I,J,K$ be a local basis of $Q$ around a point $p\in M$ as in the Definition \ref{def2}. Then the $p$-fiber of $Z$ consists precisely of those linear combinations of $I,J,K$ which are complex structures on $T_pM$. We call $Z$ the twistor space of $M$. For a point $x\in Z_p:=\pi^{-1}(p)$ we denote by $I(x)$ the complex structure defined on $T_pM$. One can prove that for a general quaternion-K\"ahler manifold $M$ the twistor space $Z$ admits a natural complex structure, however for quaternion-K\"ahler spaces with positive Ricci curvature one has even more as was shown by Salamon (\cite{Salamon}):
\begin{thrm}Let $M$ be a quaternion-K\"ahler manifold with positive Ricci curvature and $\dim M>4$. Then its twistor space $Z$ admits a K\"ahler-Einstein metric with positive Ricci curvature such that the projection $\pi:Z\to M$ is a Riemannian submersion with totally geodesic fibers.
\end{thrm}
The natural complex structure $\mathcal{J}$ of $Z$ is pretty easy to describe. Let $\mathcal{H}$ (resp. $\mathcal{V}$) denote the horizontal (resp. vertical) subbundle of the Riemannian submersion $\pi: Z\to M$. A point $z\in Z$ corresponds to a complex structure $I(z)$ on $T_{\pi(z)}M$. This gives a complex structure on $\mathcal{H}_z$ via $\pi_*$. On $\mathcal{V}_z$ we get a complex structure by the canonical complex structure of the fiber $S^2$.

Now if $M$ is a Wolf space then its twistor space is even homogeneous. A Wolf space is a quaternion-K\"ahler symmetric space of compact type and it can be written as $M=G/K\kdot \Sp(1)$, where $G$ is a simple Lie group. The twistor space of $M$ is then given by $Z=G/K\kdot \U(1)$ for some $U(1)\subset \Sp(1)$ and the twistor fibration is just the homogeneous fibration
$$\Sp(1)/\U(1)\hookrightarrow G/K\kdot \U(1)\to G/K\kdot \Sp(1).$$
For a local section $I$ of Z one can then define a local two-form $\omega_I$ by setting $\omega_I(\cdot,\cdot):=g(I\cdot,\cdot)$. Analogously to the symplectic case we define:
\begin{defin} Let $N$ be a submanifold of a quaternion-K\"ahler manifold $M$. We call $N$ a quaternion-coisotropic submanifold if and only if $$T_pN^{\perp_{\omega_I}}\subset T_pN$$ for all $p\in N$ and all sections $I$ of $Z$ locally defined around $p$. As in the symplectic case this is equivalent to saying that all normal spaces of $N$ are totally real with respect to $Z$. We call a Lie group action quaternion-coisotropic if and only if all principal orbits are quaternion-coisotropic.
\end{defin}

In the next section the twistor space will be useful for us in the analysis of K\"ahler submanifolds of a quaternion-K\"ahler manifold. The following concepts are due to Alekseevsky and Marchiafava \cite{Aleks}, \cite{Aleks2}. In the following $M$ is assumed to be a positive quaternion-K\"ahler manifold, i.~e.~the Ricci curvature of $M$ is positive.
\begin{defin} A submanifold $N$ of a quaternion-K\"ahler manifold $(M,g,Q)$ is called almost complex if there exists a section $J\in \Gamma(\left.Q\right|_N)$ such that\renewcommand{\labelenumi}{(\roman{enumi})} \begin{enumerate}\item$ J^2=-Id,$ \item $JTN=TN$.
\end{enumerate}
If $\left.J\right|_{TN}$ is integrable then $(N,J)$ is a complex manifold.
\end{defin}
\begin{defin} An almost complex submanifold $(N,J)$ of $(M,g,Q)$ is called supercomplex if
$$\nabla_{JX}J-J\nabla_XJ=0 \;\;\text{ for all } X\in TN,$$
where $\nabla$ is of course the Levi-Civita connection of $(M,g)$.
\end{defin} 
We will need the following two results:
\begin{prop}[Theorem 1.1 of \cite{Aleks}] Let $(N,J)$ be an almost complex submanifold of M with $\dim N>2$. Then $N$ is complex if and only if it is supercomplex.
\end{prop}
\begin{remark}Note that an almost complex manifold of dimension two is automatically complex.
\end{remark}
\begin{prop}[Corollary 4.6 of \cite{Aleks2}] Let $X$ be a complex submanifold of the twistor space $Z$ of $M$ and assume that the restricted twistor projection\break  $\pi:X\to Y:=\pi(X)$ is a diffeomorphism. Then $Y$ is a supercomplex manifold of $M$ and any supercomplex manifold of $M$ has such form.
\end{prop}

Let $(X,\omega)$ be a compact simply-connected irreducible homogeneous K\"ahler manifold and $H$ a compact connected Lie group acting polarly on $M$. Then it is known that the $H$-action is coisotropic (\cite{Podesta}). In other words the sections for the $H$-action are totally real with respect to the complex structure $J$, i.~e.~$T_p\Sigma\perp JT_p\Sigma$ for all sections $\Sigma\subset M$ and all points $p\in\Sigma$. 

Now let $M$ be a Wolf space endowed with a polar action of a compact connected Lie group $H$. Although there does not exist a complex structure on $M$, we can still ask a similar question:  Let $Z\subset Q\subset \End(TM)$ be the twistor bundle. Are the sections of the $H$-action totally real with respect to $Z$? More precisely we would like to prove:
\begin{thrm}\label{qc}Let $M$ be a Wolf space and $H$ a connected compact Lie group acting polarly on $M$. Then the $H$-action is quaternion-coisotropic, i.~e.~the sections are totally real with respect to $Z$.

\end{thrm}
Let us briefly recall the proof of the K\"ahler case (\cite{Podesta}). There the key observation is, that the $H$-action on $X$ admits a complex orbit. Hence there exists a polar slice  representation on a complex vector space and basically Dadoks Theorem allows to show that the sections for the slice representation are totally real. So the first idea in the quaternion-K\"ahler case would be to search for a quaternionic orbit of the $H$-action. We cannot guarantee the existence of a quaternionic orbit, however we will prove the following result:

\begin{thrm}\label{quaternionic} Let $M$ be a Wolf space 
and $H$ a compact connected Lie group acting isometrically on $M$. Then the $H$-action admits at least one of the following orbits:
\renewcommand{\labelenumi}{(\roman{enumi})}
\begin{enumerate}\item a complex orbit,
\item a totally geodesic quaternionic orbit, 
\item a $\mathbb{Z}_2$-quotient of a complex orbit in the twistor space $Z$ of $M$.
\end{enumerate}
\end{thrm}
The key idea for the proof of this result is to lift the $H$-action to a holomorphic action on the twistor space $Z\to M$. There one can find a complex orbit (cf.~the appendix) and project this orbit down to $M$ to find a good candidate for a quaternionic orbit.
\begin{lemma} Let $M$ be a Wolf space
and $\pi:Z\to M$ the twistor fibration. Let $H$ be a compact connected Lie group acting isometrically on $M$. Then the $H$-action lifts to an equivariant action on $Z$ and there exists a complex $H$-orbit in $Z$. 
\end{lemma}
\begin{proof} Let $G$ denote the identity component of the isometry group of $M$. Then an element $g\in G$ acts on $Z$ by sending a complex structure $J$ to the complex structure $\tilde{g}(J)=g_*Jg_*^{-1}$. Nitta and Takeuchi \cite{Nitta} have proven that the correspondence $g\leftrightarrow \tilde{g}$ realizes $G$ as a real form of the identity component of holomorphic automorphisms of the contact structure on $Z$ (see also \cite{Poon}). This means that $H^\mathbb{C}$ acts on $Z$ and if $H^{\mathbb{C}}=H\kdot S$ is the Iwasawa decomposition of $H^\mathbb{C}$ then $S$ has a fixed point $z$ in $Z$  by the Borel fixed point theorem. The $H$-orbit through $z$ is complex since $H^\mathbb{C}\kdot z=H\kdot S\kdot z=H\kdot s$. Actually this is only a rough sketch of the proof, since $H$ need not  be semisimple and hence the  Iwasawa decomposition is not defined. However this is only a technical problem  since one can pass to a certain cover of $H$. For the full details confer \cite{Tebege}.
\end{proof}
\begin{lemma}\label{twistor} Let $M$ be a Wolf space and $\pi: Z \to M$ the twistor fibration. Assume that a compact connected Lie group $H$ acts equivariantly with respect to $\pi$ on $M$ and $Z$. If for $z\in Z$ the orbit $H\kdot z$ is complex and $\pi^{-1}(\pi(z))\subset H\kdot z$ then all complex structures in the fiber $\pi^{-1}(\pi(z))$ leave the tangent space in $\pi(z)$ of $H\kdot\pi(z)$ invariant.
\end{lemma}
\begin{proof} Let $x\in \pi^{-1}(\pi(z))$ and $X\in T_{\pi(z)}H\kdot\pi(z)=T_{\pi(x)}H\kdot\pi(x)=\pi_*T_xH\kdot x$. Choose a vector $X'\in T_xH\kdot x$ with $X=\pi_*X'$. We then have
$$I(x)X=I(x)\pi_*X'=\pi_*\mathcal{J}X'$$ and hence $I(x)X\in T_{\pi(z)}H\kdot \pi(z)$, since $\mathcal{J}X'\in T_xH\kdot x$.
\end{proof}
We now come to the proof of Theorem \ref{quaternionic}:
\begin{proof}[Proof of Theorem \ref{quaternionic}] Let $\pi: Z \to M$ be the twistor fibration. We have seen that the $H$-action lifts to an action on $Z$ and that there exists a point $z_0\in Z$ such that $H\kdot z_0$ is complex. We set $p_0:=\pi(z_0)$ and we will prove that the orbit  $H\kdot p_0\subset M$ is of the required type. Let us first assume that the restriction $\hat{\pi} : H\kdot z_0 \to H\kdot p_0$ of $\pi$ to the orbit in $Z$ is a diffeomorphism. 
Then Corollary 4.6 of \cite{Aleks2} shows that $H\kdot p_0$ is a supercomplex submanifold of $M$ and Theorem 1.1 of \cite{Aleks} implies that $H\kdot p_0$ is a complex submanifold of $M$.

So let us now assume that the restriction of $\pi$ is not  a diffeomorphism. There are two possibilities for $\hat{\pi}$: Either there exists a point $x\in H\kdot z_0$ such that the differential of $\hat{\pi}$ in this point has not full rank or $\hat{\pi}$ is not injective. Let us first assume the former one. Without loss of generality we can assume that $x=z_0$. We will prove that in this case $H\kdot p_0$ is a quaternionic orbit. Let $v\in \ker d\hat{\pi}_{z_0}$ be a non-trivial vector. Since the orbit $H\kdot z_0$ is complex we have $\mathcal{J}v\in T_{z_0}H\kdot z_0$. Then $$d\hat{\pi}_{z_0}\mathcal{J}v=I(z_0)d\hat{\pi}_{z_0}v=0$$ shows that the vertical part $\mathcal{V}_{z_0}$ of $T_{z_0}Z$ is contained in $\ker d\hat{\pi}_{z_0}$. Since $\pi$ is equivariant with respect to the $H$-action we even have $\mathcal{V}_z\subset T_zH\kdot z_0$ for all $z\in H\kdot z_0$. This shows that for all $z\in H\kdot z_0$ the $\pi$-fiber thorugh $z$ is contained in $H\kdot z_0$. Then Lemma \ref{twistor} immediately shows that the orbit $H\kdot p_0$ is quaternionic. The fact that $H\kdot p_0$ is totally geodesic follows from the well known fact that a quaternionic submanifold is automatically totally geodesic (cf.~\cite{Gray}). Now we come to the last case: The restriction $\hat{\pi}$ is not a diffeomorphism but its differential has full rank everywhere. This means that $\hat{\pi}$ is not injective, i.~e.~there exists a point $x\in H\kdot z_0$ such that the cardinality of the set $\hat{\pi}^{-1}(\hat{\pi}(x))\cap H\kdot z_0$ is bigger than one. Assume that there exist two points $z_1,z_2\in \hat{\pi}^{-1}(\hat{\pi}(x))\cap H\kdot z_0$ such that the complex structures $I(z_1)$ and $I(z_2)$ are linearly independent. Then $I(z_1)I(z_2)$ is a third linearly independent complex structure and since all these complex structures leave $T_{p_0}H\kdot p_0$ invariant it follows that the orbit $H\kdot p_0$ is quaternionic. Clearly there exist two linearly independent complex structures in  $\hat{\pi}^{-1}(\hat{\pi}(x))\cap H\kdot z_0$ as soon as the cardinality of this set is greater or equal to three. So let us assume that the cardinality is precisely two.
Moreover we can assume that the two distinct points in $\hat{\pi}^{-1}(\hat{\pi}(x))\cap H\kdot z_0$ are antipodal to each other because otherwise the corresponding complex structures would be linearly independent. It follows that $\hat{\pi}: H\kdot z_0\to H\kdot p_0$ is a two-fold covering since $\hat{\pi}$ is a surjective local diffeomorphism and $H\kdot z_0$ is a complete connected manifold.
\end{proof}
As long as this proposition cannot be improved a direct generalization of the proof of \cite{Podesta} is not possible. However Theorem \ref{qc} can be proved if one uses the following recent classification result of Kollross:
\begin{thrm}Let $M$ be a compact symmetric space of rank greater than one whose isometry group $G$ is simple. Let $H\subset G$ be a closed connected non-trivial subgroup acting polarly on $M$. Then the action of $H$ on $M$ is hyperpolar, that is, the sections are flat in the
induced metric. Moreover, the sections are embedded submanifolds.
\end{thrm}
The higher rank case follows then from the next Lemma.
\begin{lemma} Let $H$ be a connected compact Lie group acting hyperpolarly on a Wolf space $M=G/K$. Then the sections are totally real with respect to the quaternionic structure.
\end{lemma}
\begin{proof}Let $\frak g=\frak k\oplus \frak p$ denote the corresponding Cartan decomposition and $\frak{sp}(1)\subset \frak k$ the subalgebra defining the quaternionic structure, which means there exist three generators $z_1,z_2,z_3\in \frak{sp}(1)$ such that $$\ad_{z_i}:\frak p\to \frak p \;\;\;i=1,2,3$$ are complex structures in the origin defining the quaternionic structure. Let $\Sigma \subset M$  be a section for the $H$-action and $\frak m\subset \frak p$ its tangent space in the origin. It follows $$\left<\frak m, \ad_{z_i}(\frak m)\right>=\left<\frak m, [z_i,\frak m]\right>=\left<z_i, [\frak m,\frak m]\right>=0$$ which proves the lemma.
\end{proof}
The proof of Theorem \ref{qc} follows then from the following proposition:
\begin{prop}\label{rankone} Let $M$ be a rank one quaternion-K\"ahler symmetric space, i.~e.~$M=\mathbb{P}_n(\mathbb{H})$ and $H$ a connected compact Lie group acting polarly on $M$. Then the sections of the $H$-actions are totally real with respect to the quaternionic structure.
\end{prop}
Before we come to proof we cite a result of Podest\`a and Thorbergsson \cite{Podesta2} which we will need:
\begin{thrm} Let $H$ be a connected compact Lie group acting polarly on $\mathbb{P}_n(\mathbb H)$ $(n>1)$. Then the $H$-action is orbit equivalent to an action induced by the isotropy representation of $(U,K)=\Pi_{i=1}^k(U_i,K_i)$ such that the factors $U_i/K_i$ are quaternion-K\"ahler symmetric spaces, where at least $k-1$ of them have rank one. 
\end{thrm}
\begin{proof}[Proof of Proposition \ref{rankone}] We briefly recall how the quaternionic structure on $M=\mathbb{P}_n(\mathbb{H})$ can be constructed from the one of $\mathbb{H}^{n+1}$. We write $M$ as the quotient \break $(\mathbb{H}^{n+1}\setminus \{0\})/\mathbb{H}^*$ with respect to the natural right multiplication of $\mathbb{H}^*$ on $\mathbb{H}^{n+1}\setminus\{0\}$. We can identify $\mathbb{H}^{n+1}$ with $\mathbb{R}^{4n+4}$, so let $\left<\cdot,\cdot\right>$ denote the standard inner product. We then get an $\Sp(1)$-principal bundle $$\Sp(1)\hookrightarrow S^{4n+3}\overset{\pi}{\to} \mathbb{P}_n(\mathbb{H}).$$
Let $p\in \mathbb{P}_n(\mathbb{H})$ be an arbitrary point and $z\in \pi^{-1}(p)$. The quaternionic structure of $\mathbb{P}_n(\mathbb{H})$ in the point $p$, i.~e.~on the tangent space $T_p\mathbb{P}_n(\mathbb{H})$ can be constructed as follows: Let $\mathcal{H}_z\subset T_zS^{4n+3}$ denote the horizontal space with respect to $\left<\cdot,\cdot \right>$, i.~e.~$$\mathcal{H}_z:=\ker d\pi_z^\perp=\{w\in\mathbb{H}^{n+1}\mid 0=\left<z,w\right>=\left<zi,w\right>=\left<zj,w\right>=\left<zk,w\right>\},$$ then $\mathcal{H}_z$ is invariant under right multiplication with $i,j,k$ and this gives via the isomorphism $d\pi_z: \mathcal{H}_z\to T_p\mathbb{P}_n(\mathbb{H})$ the desired quaternionic structure: For example one defines $I: T_p\mathbb{P}_n(\mathbb{H})\to T_p\mathbb{P}_n(\mathbb{H})$ by setting for $X\in T_p\mathbb{P}_n(\mathbb{H})$: $IX:=d\pi_z(\hat{X}i)$, where $\hat{X}\in \mathcal{H}_z$ with $d\pi_z\hat{X}=X$.

\noindent Consider first the case $\P_n(\mathbb{H})$ with $n>1$. Let $\Sigma\subset \mathbb{P}_n(\mathbb{H})$ be a section for the $H$-action. The above cited result of Podest\`a and Thorbergsson implies that there exist a symmetric pair $(U,K)$ such that we can assume that $\Sigma$ is the projectivization $\mathbb{P}(\frak a)$ of a linear section $\frak a\subset \frak p\cong \mathbb{H}^{n+1}$ for the $K$-action on $\frak p$, where $\frak u=\frak k\oplus\frak p$ is the corresponding Cartan decomposition. Now we decompose $\frak u$ into its irreducible components $\frak u=\sum_{l=1}^m \frak u_i$ and consider the corresponding Cartan decompositions $\frak u_l=\frak k_l\oplus \frak p_l$. The  complex structures $I,J,K: \frak p \to \frak p$ defining the quaternionic structure decompose into $I_l,J_l,K_l: \frak p_l \to \frak p_l$ such that each of them has the form $\ad_w$ for some $w\in \frak k_l$.  Let $z\in \frak a \cap S^{4n+3}\subset \frak p\cong\mathbb{H}^{n+1}$. 

\noindent We claim that $z^\perp\cap \frak a\subset \mathcal{H}_z$. Let $x\in z^\perp\cap \frak a$. We will have to prove $$0=\left<Iz,x\right>=\left<Jz,x\right>=\left<Kz,x\right>,$$ but this follows easily: Decompose $\frak a=\sum_{l=1}^m\frak a_l$ such that each $\frak a_l$ is a section for the polar action of $K_l$ on $\frak p_l$ and write $$x=\sum_{l=1}^m x_l,\;\; z=\sum_{l=1}^m z_l,\;\; I_l=\ad_{w_l} \text{ with } x_l, z_l\in \frak a_l,w_l\in \frak k_l.$$ We then get
$$\left<Iz,x\right>=\left<I\sum_{l=1}^m z_l,\sum_{r=1}^m x_r\right>=\sum_{l=1}^m\left<I_lz_l,x_l\right>=\sum_{l=1}^m\left<[w_l,z],x\right>=0.$$
The same calculations with $J$ and $K$ prove the claim.

\noindent Now let $X,Y\in T_{\pi(z)}\Sigma=T_{\pi(z)}\mathbb{P}(\frak a)$ and $X',Y' \in z^\perp\cap \frak a\subset \mathcal{H}_z$ with $d\pi_zX'=X$ and $d\pi_zY'=Y$. Let $g$ denote the Riemannian metric of $\mathbb{P}_n(\mathbb{H})$ and $\tilde{I}$ the complex structure on $T_{\pi(z)}\mathbb{P}_n(\mathbb{H})$ defined by $I$. We get 
$$g(\tilde{I}X,Y)=g(d\pi_z IX',d\pi_z Y')=\left<IX',Y'\right>=0,$$ i.~e.~$T_{\pi(z)}\Sigma$ is totally real with respect to $\tilde{I}$. Of course the argument carries over to the other two complex structures.

Now assume $n=1$, i.~e.~$M\cong S^4$. Polar actions on spheres are restrictions of polar representations. If follows from Dadoks result \cite{Dadok} that all polar actions on $S^4$ have cohomogeneity one and hence they are trivially quaternion-coisotropic.
\end{proof}


\newpage

\begin{thebibliography}{99} \addcontentsline{toc}{chapter}{Bibliography}
\bibitem[AlAl 1993]{Aleks3}D.V. Alekseevsky and A. Alekseevsky {\it Asystatic $G$-manifolds}, Differential geometry and topology (Alghere, 1992), World Sci. Publishing, River Edge, NJ, (1993), 1-22.
\bibitem[AlMa 2001]{Aleks} D.V. Alekseevsky and S. Marchiafava, {\it Almost Hermitian and K\"ahler submanifolds of a quaternionic K\"ahler manifold}, Osaka J. Math. {\bf 38} (2001) 869-904.
\bibitem[AlMa 2005]{Aleks2} D.V. Alekseevsky and S. Marchiafava, {\it A twistor construction of K\"ahler submanifolds of a quaternionic K\"ahler manifold}, Ann. Mat. Pura Appl. {\bf 184} (2005), 53-74.
\bibitem[Ak 1995]{Akhiezer}D. N. Akhiezer, {\it Lie Group Actions in Complex
Analysis,} Vieweg \& Sohn Verlagsgesellschaft, Braunschweig/Wiesbaden,
1995.
\bibitem[At 1982]{Atiyah} M. F. Atiyah, {\it Convexity and commuting Hamiltonians}, Bull. London Math. Soc. {\bf 14} (1982), 1-15.
\bibitem[Bi 2006]{Biliotti} L. Biliotti, {\it Coisotropic and polar actions on compact irreducible Hermitian symmetric spaces}, Trans. Amer. Math. Soc. {\it 358} (2006), 3003-3022.
\bibitem[BiGo 2005]{Gori} L. Biliotti and A. Gori, {\it Coisotropic and polar actions on complex Grassmannians}, Trans. Amer. Math. Soc. {\it 357} (2005), 1731-1751.
\bibitem[BaNa 1990]{Barth} W. Barth and R. Narasimhan (Eds.), {\it Several
Complex Variables VI}, Encyclopaedia of Mathematical Sciences Volume 69,
Springer Verlag, Berlin-Heidelberg-New York, 1990.
\bibitem[Be 1987]{Besse} A. L. Besse, {\it Einstein Manifolds}, vol 10 of Ergebnisse der Mathematik und Grenzgebiete, 3. Folge, Springer, Berlin, Heidelberg, New York, 1987.
\bibitem[Bo 1991]{Borel}A. Borel, {\it Linear algebraic groups}, second edition, Graduate Texts in Mathematics, No. 126, Springer-Verlag, New York, 1991. 
\bibitem[Bou 1968]{Bourbaki} Bourbaki, {\it Groups et alg{\`e}bres  de Lie,
Chap. 3}, Hermann, Paris, 1968.
\bibitem[BoSa 1958]{Bott} R. Bott and H. Samelson, {\it Applications of the theory
of Morse to symmetric spaces}, Amer. J. Math {\bf 80} (1958), 964-1029. {\it
Correction} in Amer. J. Math. {\bf 83} (1961), 207-208.
\bibitem[Br 1972]{Bredon2} G. E. Bredon, {\it Introduction to Compact
Transformation Groups}, Academic Press, New York and London, 1972. 
\bibitem[Co 1971]{Conlon} L. Conlon, {\it Variational completeness and $K$-transversal domains}, J. Differential Geom. {\bf 5} (1971), 135-147.
\bibitem[Dad 1985]{Dadok} J. Dadok, {\it Polar coordinates induced by actions of compact Lie groups}, Trans. Amer. Math. Soc. {\bf 288} (1985), 123-137.
\bibitem[Da 1996]{Da} R. Dabrowski, {\it On normality of the closure of a generic torus orbit in $G/P$}, Pacific J. of Math {\bf 172} (1996), 321-330.
\bibitem[DuKo 2000]{Duistermaat} J. J. Duistermaat and J. A. C. Kolk, {\it Lie
Groups}, Springer Verlag, New York-Heidelberg-Berlin, 2000. 
\bibitem[FlHa 1991]{FlHa} H. Flaschka and L. Haine, {\it Torus orbits in $G/P$}, Pacific. J. Math., {\bf 149} (1991), 251-292.
\bibitem[GeSe 1987]{GeSe} I.M. Gel'fand and V.V. Serganova, {\it Combinatorial geometries and torus strata on homogeneous compact manifolds}, Russian Math. surveys {\bf 42} (1987), 133-168.
\bibitem[GoPo 2006]{Gori2} A. Gori and F. Podest\`a, {\it Symplectically asystatic actions of compact Lie groups}, Transform. Groups {\bf 11} (2006), 177-184.
\bibitem[Go 2004]{Go} C. Gorodski, {\it Polar actions on compact symmetric spaces
which admit a totally geodesic principal orbit}, Geometriae Dedicata  {\bf 103} (2004), 193-204.
\bibitem[Gr 1965]{Gray} A. Gray, {\it A note on manifolds whose holonomy is a subgroup of \break$\Sp(n)\kdot\Sp(1)$}, Mich. Math. J. {\bf 16} (1965).
\bibitem[GuSt 1982]{Guillemin} V. Guillemin and S. Sternberg, {\it Convexity properties of the moment mapping}, Invent. Math. {\bf 67} (1982), 491-513.
\bibitem[GuSt 1984]{Guillemin 2} V. Guillemin and S. Sternberg, {\it Symplectic Techniques in Physics}, Cambridge University Press, Cambdrige, 1984.
\bibitem[GuRo 1965]{Gunning} R. C. Gunning and H. Rossi, {\it Analytic Functions
in Several Complex Variables}, Prentice-Hall, Englewood cliffs, N.J., 1965.
\bibitem[Ha 1977]{Hartshorne} R. Hartshorne, {\it Algebraic Geometry}, Springer
Verlag, New York-Heidelberg-Berlin, 1977.
\bibitem[HPTT 1994]{Heintze} E. Heintze, R. S. Palais, C.-L. Terng and G.
Thorbergsson, {\it Hyperpolar actions and k-flat homogeneous spaces}, J. reine
angew. Math. {\bf 454} (1994), 163-179.
\bibitem[He 1962]{Helgason} S. Helgason, {\it Differential Geometry and
Symmetric Spaces}, Academic Press, New York and London, 1962.
\bibitem[Ho 1965]{Hochschild} G. Hochschild, {\it The Structure of Lie
Groups}, Holden-Day, San Francisco, London, Amsterdam, 1965.
\bibitem[HuWu 1990]{HuWu} A. T. Huckleberry and T. Wurzbacher, {\it Multiplicity-free complex manifolds} Math. Annalen {\bf 286} (1990) 261-280.
\bibitem[Hu 1975]{Humphreys} J. E. Humphreys, {\it Linear algebraic groups}, Graduate Texts in Mathematics, No. 21, Springer-Verlag, New York-Heidelberg, 1975.
\bibitem[Ki 1984]{Kirwan} F. Kirwan, {\it Convexity properties of the moment mapping III}, Invent. Math. {\bf 77} (1984), 547-552.
\bibitem[Kn 1996]{Knapp} A. W. Knapp, {\it Lie Groups Beyond an Introduction},
Birkh{\"a}user, Boston, Basel, Berlin, 1996	
\bibitem[Ko 1972]{Kobayashi} S. Kobayashi, {\it Transformation Groups in
Differential Geometry}, Springer Verlag, Berlin-Heidelberg-New York, 1972.
\bibitem[KoNo 1963]{KN} S. Kobayashi and K. Nomizu, {\it Foundations of
Differential Geometry}, vol. I, Interscience Publishers, J. Wiley \& Sons,
1963.
\bibitem[KoNo 1969]{KN2} S. Kobayashi and K. Nomizu, {\it Foundations of
Differential Geometry}, vol. II, Interscience Publishers, J. Wiley \& Sons,
1969.
\bibitem[Kol 1998]{Kollross} A. Kollross, {\it A classification of hyperpolar and
cohomogeneity one actions}, PhD thesis, Augsburg (1998).
\bibitem[Kol 2006]{Kollross2} A. Kollross, {\it Polar actions on symmetric spaces}, arXiv:math.DG/0506312, to appear in Journal of Differential Geometry.
\bibitem[Kos 1955]{Kostant} B. Kostant, {\it Holonomy and the Lie algebra in
infinitesimal motions of a Riemannian manifold}, Trans. Amer. Math. Soc.
{\bf 80} (1955), 528-542.
\bibitem[MSZ 1956]{Montgomery1} D. Montgomery, H. Samelson and L. Zippin, {\it
Singular points of a compact transformation group}, Ann. of Math. {\bf 63}
(1956), 1-9.
\bibitem[MSY 1956]{Montgomery2} D. Montgomery, H. Samelson and C. T. Yang, {\it
Exceptional orbits of highest dimension}, Ann. of Math. {\bf 64}
(1956), 131-141.
\bibitem[NiTa 1987]{Nitta} T. Nitta and M. Takeuchi, {\it Contact structures on twistor spaces}, J. Math. Soc. Japan {\bf 39} (1987) 139-162.
\bibitem[Oda 1988]{Oda} T. Oda, {\it Convex bodies and algebraic geometry - an introduction to the theory of toric varieties}, Springer-Verlag, Berlin-Heidelberg-New York, 1988.
\bibitem[Oda 1991]{Oda2}T. Oda, {\it Geometry of toric varieties}, Proc. of the Hyderabad conference on algebraic groups, Manoj Prakashan, Madra-India, 1991.
\bibitem[PaTe 1987]{Palais1} R. S. Palais and C.-L. Terng, {\it A general theory of
canonical forms}, Trans. Amer. Math. Soc. {\bf 300} (1987), 771-789.
\bibitem[PaTe 1988]{Palais2} R. S. Palais and C.-L. Terng, {\it Critical point
theory and submanifold geometry}, Lecture Notes in Mathematics 1353,
Springer-Verlag, Berlin-NewYork, 1988.
\bibitem[PoTh 1999]{Podesta2} F. Podest{\`a} and G. Thorbergsson, {\it Polar
actions on rank one symmetric spaces,} J. Differential Geom. {\bf 53} (1999),
131-175.
\bibitem[PoTh 2002]{Podesta} F. Podest{\`a} and G. Thorbergsson, {\it Polar and
coisotropic actions on K{\"a}hler manifolds}, Trans. Amer. Math. Soc. {\bf
354} (2002), 1759-1781.
\bibitem[PoTh 2003]{Podesta3} F. Podest\`a and T. Thorbergsson, {\it Coisotropic actions on compact homogeneous K\"ahler manifolds}. Math. Z. {\bf 243} (2003), 471-490.
\bibitem[PoSa 1991]{Poon} Y.S. Poon and M. Salamon, {\it Quaternionic K\"ahler 8-manifolds with positive scalar curvature}, J. Differential Geom. {\bf 33} (1991), 363-378.
\bibitem[Sa 1992]{Sakai} T. Sakai, {\it Riemannian Geometry}, AMS Providence,
Rhode Island, 1992.
\bibitem[Sa 1982]{Salamon} S. M. Salamon, {\it Quaternionic K\"ahler manifolds}, Invent. Math {\bf 67} (1982), 143-171.
\bibitem[Sa 1999]{Salamon2}S. M. Salamon, {\it Quaternion K\"ahler geometry}, Surv. Differ. Geom. VI, Int. Press, Boston, MA, (1999), 83-121
\bibitem[SCM 1997]{Sanchez} C. U. S\'anchez, A. L. Cal\'i and J. L. Moreschi, {\it Spheres in Hermitian symmetric spaces and flag manifolds}, Geom. Dedicata {\bf 64} (1997), 261-276. 
\bibitem[Te 2006]{Tebege} S. Tebege, {\it Polar actions on Hermitian and quaternion-K\"ahler symmetric spaces}. PhD thesis. K\"oln (2006).
\bibitem[Va 1984]{Varadarajan} V.S. Varadarajan, {\it Lie Groups, Lie Algebras,
and Their Representations}, Springer Verlag, New
York-Berlin-Heidelberg-Tokyo, 1984.
\bibitem[Wo 1965]{Wolf} J. A. Wolf, {\it Complex homogeneous contact manifolds and quaternionic symmetric spaces}, J. Math. Mech., {\bf 14} (1965), 1033-1047.
\end{thebibliography}
\end{document}